\documentclass[preprint,12pt,nopreprintline]{elsarticle}

\usepackage{amsmath,amssymb,amsthm}
\usepackage{graphicx}
\usepackage{booktabs}
\usepackage{microtype}
\usepackage{mathtools}
\usepackage{enumitem}
\usepackage{bm}
\usepackage{placeins}
\usepackage{url}
\usepackage{hyperref}

\biboptions{authoryear}

\newtheorem{theorem}{Theorem}
\newtheorem{lemma}{Lemma}

\DeclareMathOperator{\Var}{Var}

\DeclareMathOperator{\adj}{adj}

\newcommand{\E}{\mathbb{E}}
\newcommand{\R}{\mathbb{R}}

\begin{document}
\begin{frontmatter}

\title{Geary-type expectation identities for Mardia's sample
skewness and kurtosis}

\author[ous]{Shigekazu Nakagawa\corref{cor1}}
\ead{s-nakagawa@ous.ac.jp}
\author[tus]{Hiroki Hashiguchi}
\ead{hiro@rs.tus.ac.jp}
\author[otsuma]{Yoko Ono}
\ead{onoyk@otsuma.ac.jp}
\cortext[cor1]{Corresponding author.}

\address[ous]{Institute for the Advancement of Higher Education,
Okayama University of Science, 1-1 Ridai-cho, Okayama 700-0005, Japan}
\address[tus]{Faculty of Science Division I, Tokyo University of Science,
1-3 Kagurazaka, Shinjuku-ku, Tokyo 162-8601, Japan}
\address[otsuma]{Faculty of Data Science, Otsuma Women's University,
12 Sanbancho, Chiyoda-ku, Tokyo 102-8357, Japan}

\begin{abstract}
Under multivariate normality, we establish exact Geary-type expectation
identities for Mardia's sample skewness and kurtosis.
A polar decomposition of the centered Gaussian data matrix shows that
each of these statistics is independent of the sample covariance matrix.
This independence allows the moments of each statistic to be expressed
as ratios of expectations involving polynomial numerators and powers
of the determinant of the sample covariance matrix.
The identities hold in arbitrary dimension $p$.
To illustrate their use, we give exact expressions for the variance of Mardia's sample skewness 
and the third central moment of Mardia's sample kurtosis for \(p=2\) and \(p=3\).
As an application of the exact variance formulas, we construct a moment-matched chi-square approximation 
to the null distribution of Mardia's sample skewness 
and compare its performance with Mardia's classical and mean-corrected approximations.
\end{abstract}


\end{frontmatter}

\section{Introduction}\label{sec:introduction}

Mardia's measures of sample skewness and kurtosis are widely used to
assess multivariate normality \citep{Mardia1970}. Their sampling
distributions are complicated, and higher-order moments are important
for obtaining accurate distributional approximations. Although exact
expressions for several lower-order moments are available
\citep{Mardia1974}, direct higher-order moment calculations are difficult
because both statistics involve the inverse of a random sample
covariance matrix.
\citet{MardiaKanazawa1983} derived the leading term
in an asymptotic expansion of the third central
moment of Mardia's sample kurtosis,
while \citet{MardiaFoster1983} investigated the
correlation between Mardia's sample skewness
and kurtosis.

In the univariate normal setting, \citet{Geary1933} established identities
that express expectations of ratios
arising in moments of sample skewness
and kurtosis as ratios of expectations.
Specifically, let
$X_1,\ldots,X_n$ be a univariate normal sample, and write
$m_k=n^{-1}\sum_{i=1}^n(X_i-\bar X)^k$ for $k=2,3,4$. With
$b_{1,1}=m_3^2/m_2^3$ and $b_{2,1}=m_4/m_2^2$, these identities are
\begin{equation}\label{eq:geary_univariate}
 \E[b_{1,1}^{r}]
    =\frac{\E[m_3^{2r}]}{\E[m_2^{3r}]},\qquad
 \E[b_{2,1}^{r}]
    =\frac{\E[m_4^{r}]}{\E[m_2^{2r}]},
 \quad r=1,2,\ldots.
\end{equation}
The identities avoid the need to evaluate expectations of ratios directly.
This raises the question of whether analogous identities hold for Mardia's sample skewness and kurtosis.

The main purpose of this paper is to give a rigorous $p$-variate extension
of the expectation identities under normal sampling.
The key step is to show, using the polar decomposition of a centered
Gaussian matrix, that Mardia's sample skewness and kurtosis are each
independent of the sample covariance matrix.
This independence allows the moments of each statistic to be expressed
as ratios of expectations of polynomial numerators to moments of
the determinant of the sample covariance matrix.
The denominator expectations can be evaluated in closed form using
the determinant moment formula for the Wishart distribution.

To illustrate the utility of these identities, we obtain exact expressions for the variance of 
Mardia's sample skewness and the third central moment of Mardia's sample kurtosis for \(p=2\) and \(p=3\).
The exact variance formulas, together with the known exact mean,
are then used to construct a moment-matched chi-square approximation
to the null distribution of Mardia's sample skewness.
Monte Carlo simulations compare its finite-sample null rejection probabilities with those of 
Mardia's classical and mean-corrected chi-square approximations.

\section{Independence and Geary-type expectation identities}
\label{sec:identities}

Let $\bm X_1,\ldots, \bm X_n\overset{\mathrm{iid}}{\sim} N_p( \bm \mu,\Sigma)$
with mean vector $\bm \mu$ and covariance matrix $\Sigma$, where $n>p$ and $\Sigma$ is positive definite. 
Define the sample mean and covariance matrix by
$\overline{\bm X} = (1/n)\sum_{i=1}^n\bm X_i$, 
$S = (1/n) \sum_{i=1}^n(\bm X_i-\overline{\bm X})(\bm X_i-\overline{\bm X})^\top$.
\cite{Mardia1970} defines 
the sample skewness and kurtosis by
\begin{align}
 b_{1,p}&=\frac1{n^2}\sum_{i=1}^n\sum_{j=1}^n
 \left\{( \bm X_i-\overline{ \bm X})^{\top}S^{-1}(\bm X_j-\overline{ \bm X})\right\}^{3},
 \label{eq:mardia_skewness}\\
 b_{2,p}&=\frac1n\sum_{i=1}^n
 \left\{(\bm X_i-\overline{\bm X})^{\top}S^{-1}( \bm X_i-\overline{\bm X})\right\}^{2},
 \label{eq:mardia_kurtosis}
\end{align}
respectively. We define the Stiefel manifold $V_p(\R^m)$ by 
$V_p(\R^m)=\{U\in\R^{m\times p}:U^\top U=I_p\}$ and write 
the $p$-dimensional Wishart
distribution with $m$ degrees of freedom and scale matrix $\Sigma$ as $W_p(m,\Sigma)$.

 These statistics $b_{1,p}$ and $b_{2,p}$ remain unchanged under 
 the affine transformation $\bm X_i=\Sigma^{1/2} \bm Z_i + \bm \mu$ ($i=1,\ldots,n$), 
 where $\bm Z_i \sim N_{p}(\bm 0, I_{p})$.
To prove the following independence results, it suffices to consider \(\boldsymbol\mu=\boldsymbol0\) and \(\Sigma=I_p\).

\begin{lemma}[Gaussian polar decomposition; \citealp{Chikuse2003}]
\label{lem:polar}
Let $Z$ be an $m\times p$ matrix with independent standard normal
entries, where $m\ge p$. Define
\begin{equation}\label{eq:polar}
 T=Z^\top Z,\qquad U=ZT^{-1/2},\qquad Z=UT^{1/2}.
\end{equation}
Then $U$ and $T$ are independent; $U$ is uniformly distributed on
$V_p(\R^m)$; and $T\sim W_p(m,I_p)$.
\end{lemma}
The following two theorems follow from Lemma~\ref{lem:polar}.

\begin{theorem}\label{thm:b1pindependence}
Under the preceding assumptions, 
$b_{1,p}$ is independent of $S$,
and hence of $D=\det(S)$.
\end{theorem}

\begin{theorem}\label{thm:b2pindependence}
Under the preceding assumptions, 
$b_{2,p}$ is independent of $S$,
and hence of $D=\det(S)$.
\end{theorem}

\begin{proof}[Proof of Theorems~\ref{thm:b1pindependence} and \ref{thm:b2pindependence}]
Without loss of generality, we prove the result for
$\bm{\mu}=\bm{0}$ and $\Sigma=I_p$.
Let $X \in\R^{n\times p}$ be the data matrix
whose $i$th row is $\bm X_i^\top$, and set
\begin{equation}\label{eq:projection}
 P=I_n-\frac1n\boldsymbol{1}_n\boldsymbol{1}_n^\top.
\end{equation}
Choose an $n\times(n-1)$ matrix $H$ whose columns form an
orthonormal basis of the subspace orthogonal to $\boldsymbol{1}_n$.
Then $HH^\top=P$ and $H^\top H=I_{n-1}$. 
Writing $\bm h_i=H^\top \bm e_i$,
where $\bm e_i$ is the $i$th standard basis vector of $\R^n$, we have
\begin{equation}\label{eq:reduction}
 Y=H^\top X\in\R^{(n-1)\times p},\qquad
 Y^\top Y=X^\top PX=nS.
\end{equation}
Because the columns of $X$ are independent $N_n(\bm 0, I_n)$ vectors,
the entries of $Y$ are independent $N(0,1)$ variables. In particular,
\begin{equation}\label{eq:wishart_A}
 A=Y^\top Y=nS\sim W_p(n-1,I_p).
\end{equation}
Since $n>p$, $A$ is positive definite almost surely. 
Applying Lemma~\ref{lem:polar} to $Y$, let $U=YA^{-1/2}$.
Then $Y=UA^{1/2}$, $U$ and $A$ are independent,
and $U$ is uniformly distributed on $V_p(\R^{n-1})$.

For $1\le i,j\le n$, define
$d_{ij}=(\bm X_i-\overline{\bm X})^\top S^{-1}(\bm X_j-\overline{\bm X})$. Since
$(\bm X_i-\overline{\bm X})^\top=\bm e_i^\top PX= \bm h_i^\top Y$ and $S=A/n$,
\begin{equation}
\label{eq:dij_U}
 d_{ij} = \bm h_i^\top YS^{-1}Y^\top \bm h_j
   =n\, \bm h_i^\top U A^{1/2}A^{-1}A^{1/2}U^\top \bm h_j
  =n\, \bm h_i^\top UU^\top \bm h_j.
\end{equation}
Consequently, $b_{1,p}$ and $b_{2,p}$ are given by
\begin{equation}\label{eq:stats_U}
 b_{1,p}=\frac1{n^2}\sum_{i=1}^n\sum_{j=1}^n
 (n \bm h_i^\top UU^\top \bm h_j)^3, 
 \qquad
 b_{2,p}=\frac1n\sum_{i=1}^n
 (n \bm h_i^\top UU^\top \bm h_i)^2,
\end{equation}
both of which are functions of $U$. Since $S=A/n$ is a function of $A$, 
$b_{1,p}$ is independent of $S$, and hence of $D$.
Theorem~\ref{thm:b2pindependence} 
follows similarly from the second expression in \eqref{eq:stats_U}.
\end{proof}

The independence also follows from Basu's theorem
\citep{Basu1955}, since $(\overline{\bm X},S)$ is a
complete sufficient statistic for $(\bm\mu,\Sigma)$
and each of \(b_{1,p}\) and \(b_{2,p}\) is ancillary by affine invariance.

Let \(D=\det(S)\), 
and let \(\operatorname{adj}(S)\) 
denote the adjugate of \(S\). 
Then \(S^{-1}=\operatorname{adj}(S)/D\).
Define the polynomial numerators
\begin{align}
 N_{1,p}&=\frac1{n^2}\sum_{i=1}^n\sum_{j=1}^n
 \bigl\{(\bm X_i-\overline{\bm X})^\top\adj(S)( \bm X_j-\overline{\bm X})\bigr\}^3,
 \label{eq:N1}\\
 N_{2,p}&=\frac1n\sum_{i=1}^n
 \bigl\{(\bm X_i-\overline{\bm X})^\top\adj(S)(\bm X_i-\overline{\bm X})\bigr\}^2.
 \label{eq:N2}
\end{align}
We have $b_{1,p} = N_{1,p}/D^3$ and 
$b_{2,p} = N_{2,p}/D^2$.

\begin{theorem}[Geary-type expectation identities]
\label{thm:geary}
For all nonnegative integers $r$,
\begin{equation}\label{eq:geary_p}
 \E[b_{1,p}^r]=\frac{\E[N_{1,p}^r]}{\E[D^{3r}]},\qquad
 \E[b_{2,p}^r]=\frac{\E[N_{2,p}^r]}{\E[D^{2r}]}.
\end{equation}
\end{theorem}

\begin{proof}
Theorems~\ref{thm:b1pindependence} and \ref{thm:b2pindependence} yield
$
 \E[N_{1,p}^r]
 =\E[b_{1,p}^r D^{3r}]
 =\E[b_{1,p}^r] \E[D^{3r}]
$ and
$\E[N_{2,p}^r]
 =\E[b_{2,p}^r D^{2r}]
 =\E[b_{2,p}^r] \E[D^{2r}],$
which proves \eqref{eq:geary_p}. 
\end{proof}

When \(p=1\), we have \(D=m_2\), \(N_{1,1}=m_3^2\), and \(N_{2,1}=m_4\). 
Thus, \eqref{eq:geary_p} reduces to the univariate Geary-type identities in \eqref{eq:geary_univariate}.

\citet{Mardia1974} obtained exact moments by evaluating
moments of certain bilinear forms.
Here the adjugate representation 
reduces the calculation of \(\E[b_{1,p}^r]\) to the evaluation of 
the polynomial expectation \(\E[N_{1,p}^r]\).
Since $nS\sim W_p(n-1,\Sigma)$, 
the denominator expectations are given by the Wishart determinant moment formula
\citep[Theorem~3.2.15]{Muirhead1982}: 
\begin{equation}\label{eq:det_moment}
 \E[D^q]
 =|\Sigma|^q\left(\frac2n\right)^{pq}
 \frac{\Gamma_p\bigl(\frac{n-1}{2}+q\bigr)}
      {\Gamma_p\bigl(\frac{n-1}{2}\bigr)},\qquad q\ge0,
\end{equation}
where $q$ is an integer and
$$\Gamma_p(a)=\pi^{p(p-1)/4}\prod_{j=1}^p
 \Gamma\left(a-\frac{j-1}{2}\right), \
 \Re(a)>\frac{p-1}{2}.$$
By affine invariance, we may set \(\Sigma=I_p\) when evaluating the ratios in \eqref{eq:geary_p}.

\section{Illustrative examples}\label{sec:moments}

\subsection{Exact moments for $p=2$ and $p=3$}
For $p=2$, let $m_{rs} = (1/n) \sum_{i=1}^n \left( X_{i1}-\overline X_1 \right)^r \left( X_{i2}-\overline X_2 \right)^s$,
and put $a=m_{02}$, $b=-m_{11}$, and $c=m_{20}$.
The polynomial numerators defined above admit the following
matrix representations.
With $\bm m_3=(m_{30},m_{21},m_{12},m_{03})^\top$,
\begin{equation}\label{eq:N12-matrix}
 N_{1,2}=\bm m_3^\top C\bm m_3,
 \quad
 C=\begin{pmatrix}
 a^3 & 3a^2b & 3ab^2 & b^3\\
 3a^2b & 3a(ac+2b^2) & 3b(2ac+b^2) & 3b^2c\\
 3ab^2 & 3b(2ac+b^2) & 3c(ac+2b^2) & 3bc^2\\
 b^3 & 3b^2c & 3bc^2 & c^3
 \end{pmatrix}.
\end{equation}
Moreover, with $\bm a=(a,b,c)^\top$,
\begin{equation}\label{eq:N22-matrix}
 N_{2,2}=\bm a^\top M_4\bm a,
 \quad
 M_4=\begin{pmatrix}
 m_{40} & 2m_{31} & m_{22}\\
 2m_{31} & 4m_{22} & 2m_{13}\\
 m_{22} & 2m_{13} & m_{04}
 \end{pmatrix}.
\end{equation}
By Theorem~\ref{thm:geary}, the calculation of $\E[b_{1,2}^2]$ and
$\E[b_{2,2}^3]$ reduces to evaluating $\E[N_{1,2}^2]$
and $\E[N_{2,2}^3]$, respectively.

The expectations of the polynomial numerators
were evaluated using the change-of-basis algorithm
for multisymmetric polynomials described by
\citet{NakagawaNikiHashiguchi1998}.
In Mathematica, we used the built-in function
\texttt{MomentConvert} to convert linear
combinations of products of power-sum
multisymmetric polynomials into linear
combinations of augmented multisymmetric
polynomials.
These two types of polynomials are represented
by the built-in Mathematica functions
\texttt{PowerSymmetricPolynomial} and
\texttt{AugmentedSymmetricPolynomial},
respectively.

This computation yields the following
exact expression for the second raw moment
of the bivariate sample skewness:
\begin{equation}
 \E[b_{1,2}^2]
 = \frac{288(3n^4+71n^3+5n^2-215n+40)}
 {(n+1)(n+2)(n+3)(n+5)(n+7)(n+9)}.
\label{eq:b1second}
\end{equation}
Similarly, the third raw moment of the bivariate sample kurtosis is
\begin{align}
 \E[b_{2,2}^3]
 &=\frac{512(n^6+24n^5+170n^4-96n^3-567n^2+552n+108)}
 {(n+1)(n+2)(n+3)(n+5)(n+7)(n+9)}.                                \label{eq:b2third}
\end{align}
\citet{Mardia1970, Mardia1974} obtained the following
exact moment formulas:
\begin{align}
 \E[b_{1,p}]&=\frac{p(p+2)\left\{(n+1)(p+1)-6\right\}}{(n+1)(n+3)}, \label{eq:mean_b1}\\
 \E[b_{2,p}]&=\frac{p(p+2)(n-1)}{n+1}, \\
 \Var[b_{2,p}] & =
 \frac{8p(p+2)(n-3)(n-p-1)(n-p+1)}{(n+1)^2(n+3)(n+5)}. \label{eq:var_b2}
\end{align}
Let $\mu_3[b_{2,p}]
= \E[(b_{2,p}-\E[b_{2,p}])^3]$ denote the third central moment of $b_{2,p}$.
Combining \eqref{eq:b1second} and \eqref{eq:b2third} with \eqref{eq:mean_b1}--\eqref{eq:var_b2}, we obtain the 
following exact expressions for $\Var[b_{1,2}]$ and $\mu_3[b_{2,2}]$:
\begin{align}
 \Var[b_{1,2}] = & \,
 \frac{288(n-3)(n^5+44n^4+150n^3-40n^2-151n+380)}
 {(n+1)^2(n+2)(n+3)^2(n+5)(n+7)(n+9)},
\label{eq:var_b12}
\\
\mu_{3}[b_{2,2}] 
 = & \, \frac{1024(n-3)(n^2-5n+2)(5n^3-21n^2-5n+117)}
 {(n+1)^3(n+2)(n+3)(n+5)(n+7)(n+9)}.
\label{eq:mu3_b2}
\end{align}
\cite{MardiaKanazawa1983} derived the leading term in the asymptotic expansion for $\mu_{3}[b_{2,p}]$ as
\begin{equation}
\mu_{3}[b_{2,p}] = \frac{64p(p+2)(p+8)}{n^2} + O(n^{-3}).
\label{eq:mardiakanazawa}
\end{equation}
For \(p=2\), the leading coefficient in  \eqref{eq:mardiakanazawa} is \(64\cdot2\cdot4\cdot10=5120\), 
which agrees with \(1024\cdot5=5120\) from \eqref{eq:mu3_b2}.

Similarly, for $p=3$, we obtain 
\begin{align}
\E[b_{1,3}^2] = &
\frac{540 \left(8 n^5+166 n^4+163 n^3-440 n^2-171
   n+658\right)}{(n-1) (n+1) (n+2) (n+3) (n+5) (n+7) (n+9)},
   \label{eq:b1second_p3}
\\
\E[b_{2,3}^3] = &
\frac{15 P(n)}{(n-1)(n+1)(n+2)(n+3)(n+5)(n+7)(n+9)}, 
\label{eq:b2third_p3}
\end{align}
where 
\begin{align*}
P(n) = & 225n^7+4860n^6+28379n^5-28894n^4-133489n^3\\
&+160960n^2+66101n-147294.
\end{align*}
Then, we obtain 
\begin{align}
\Var[b_{1,3}] = &
\frac{360 (n-4) (n-2) \left(2 n^5+99 n^4+452 n^3+390 n^2+26
   n+567\right)}{(n-1) (n+1)^2 (n+2) (n+3)^2 (n+5) (n+7) (n+9)},
\label{eq:var_b1_p3}
\\
\mu_{3}[b_{2,3}] = & 
\frac{960 (n-4) (n-2) \left(n^2-5 n+2\right) \left(11 n^3-68 n^2+65 n+504\right)}{(n-1) (n+1)^3 (n+2) (n+3) (n+5) (n+7) (n+9)}.
\label{eq:mu3_b2_p3}
\end{align}
For \(p=3\), the coefficient of \(n^{-2}\) obtained 
by expanding \eqref{eq:mu3_b2_p3} agrees with that in \eqref{eq:mardiakanazawa}, since \(64\cdot3\cdot5\cdot11=960\cdot11=10560\).

\subsection{Numerical verification}\label{sec:simulation}

To assess the exact formulas \eqref{eq:var_b12}, \eqref{eq:mu3_b2}, \eqref{eq:var_b1_p3}, and \eqref{eq:mu3_b2_p3} numerically, 
we generated $10^6$ independent samples of size $n$ from $N_p( \bm 0, I_p)$ for each
$n\in\{20,30,50,100\}$ and
$p\in\{2,3\}$.
Table~\ref{tab:moment_verification} compares the exact values of
$\Var[b_{1,2}]$ and
$\mu_3[b_{2,2}]$ with their Monte Carlo
estimates, indicated by hats.
Table~\ref{tab:mardia_p3_mc} presents the corresponding
results for $p=3$.
The Monte Carlo estimates are in good
agreement with the exact formulas,
allowing for simulation variability.

\begin{table}[htbp]
\centering
\caption{Exact variance and third central moment for $p=2$ and their Monte Carlo estimates (based on $10^6$ replications).}
\label{tab:moment_verification}
\begin{tabular}{rcccc}
\toprule
$n$ & $\Var[b_{1,2}]$ & $\widehat{\Var[b_{1,2}]}$ & $\mu_3[b_{2,2}]$ & $\widehat{\mu_3[b_{2,2}]}$ \\
\hline
20 & 0.5566 & 0.5549 & 1.8120 & 1.8033 \\
30 & 0.2940 & 0.2954 & 1.5188 & 1.5325 \\
50 & 0.1168 & 0.1165 & 0.9174 & 0.9137 \\
100 & 0.0301 & 0.0301 & 0.3410 & 0.3441 \\
\bottomrule
\end{tabular}
\end{table}

\begin{table}[htbp]
\centering
\caption{Exact variance and third central moment for $p=3$ and their
Monte Carlo estimates (based on $10^6$ replications).}
\label{tab:mardia_p3_mc}
\begin{tabular}{rcccc}
\toprule
$n$ & $\Var[b_{1,3}]$ & $\widehat{\Var[b_{1,3}]}$ & $\mu_3[b_{2,3}]$ & $\widehat{\mu_3[b_{2,3}]}$ \\
\hline
20 & 1.4129 & 1.4137 & 2.9992 & 3.0030 \\
30 & 0.7552 & 0.7577 & 2.7176 & 2.7390 \\
50 & 0.3004 & 0.2997 & 1.7414 & 1.7205 \\
100 & 0.0768 & 0.0766 & 0.6753 & 0.6721 \\
\bottomrule
\end{tabular}
\end{table}

\section{A moment-matched chi-square approximation}
\label{sec:chi_squared}

The exact variance formulas \eqref{eq:var_b12} and \eqref{eq:var_b1_p3}
can be used to construct a moment-matched
chi-square approximation to the null
distribution of Mardia's sample skewness.
The basic idea is to approximate the distribution of $b_{1,p}$
by a scaled chi-square distribution whose mean and variance
agree with the corresponding exact moments of $b_{1,p}$.

Let $\mu_{n,p}= \E[b_{1,p}]$ and $\sigma_{n,p}^2 = \Var[b_{1,p}]$. Suppose that $\sigma_{n,p}^2>0$.
We approximate the distribution of $b_{1,p}$ by that of
$c_{n,p}\chi^2_{\nu_{n,p}}$, where $c_{n,p}>0$ is a scale
parameter and $\nu_{n,p}>0$ is the degrees-of-freedom parameter.
Since the mean and variance of this scaled chi-square
distribution are $c_{n,p}\nu_{n,p}$ and
$2c_{n,p}^2\nu_{n,p}$, respectively, matching these
moments with $\mu_{n,p}$ and $\sigma_{n,p}^2$ gives
\begin{equation}
c_{n,p}
=
\frac{\sigma_{n,p}^2}{2\mu_{n,p}},
\qquad
\nu_{n,p}
=
\frac{2\mu_{n,p}^2}{\sigma_{n,p}^2}.
\label{eq:mm_parameters}
\end{equation}
Thus, the resulting approximation is
\begin{equation}
\frac{b_{1,p}}{c_{n,p}}
\approx \chi^2_{\nu_{n,p}}.
\label{eq:mm_chisq}
\end{equation}

For $p=2$ and $p=3$, the exact mean $\mu_{n,p}$ is
given by \eqref{eq:mean_b1}, while the exact variances are given
by \eqref{eq:var_b12} and \eqref{eq:var_b1_p3}, respectively. Hence, the
moment-matched approximation is explicitly available
in these two cases for $n>p+1$.

For comparison, \citet{Mardia1970} established the classical
asymptotic result
\begin{equation}
\frac{nb_{1,p}}{6}
\xrightarrow{d} \chi_f^2
\quad (n\to\infty),
\label{eq:mardia_classical}
\end{equation}
where $f=p(p+1)(p+2)/6$.
\citet{Mardia1974} also proposed a finite-sample
mean correction, which gives the approximation
\begin{equation}
\frac{nK_{n,p}}{6}b_{1,p}
\approx
\chi_f^2, 
\qquad
K_{n,p}
=
\frac{(p+1)(n+1)(n+3)}
{n\{(n+1)(p+1)-6\}}.
\label{eq:mardia_corrected}
\end{equation}
Using \eqref{eq:mean_b1}, it follows that
$\E\left[nK_{n,p}b_{1,p}/6\right]=f$.
Thus, the mean-corrected approximation matches
the exact mean while retaining the asymptotic
degrees of freedom $f$.
In contrast, the proposed approximation
\eqref{eq:mm_chisq} matches both the exact mean
and the exact variance by adjusting the scale
and the degrees-of-freedom parameter.

Furthermore, for $p=2$ and $p=3$, the exact
moment formulas \eqref{eq:mean_b1}, \eqref{eq:var_b12}, and \eqref{eq:var_b1_p3} yield
\[
c_{n,p}
=
\frac{6}{n}+O(n^{-2}),
\qquad
\nu_{n,p}
=
f+O(n^{-1}).
\]
Consequently, the moment-matched approximation
is consistent with Mardia's classical
chi-square limit in \eqref{eq:mardia_classical}.

\begin{figure}[tbh]
\centering
\begin{minipage}[t]{0.48\textwidth}
  \centering
  \includegraphics[width=\linewidth]{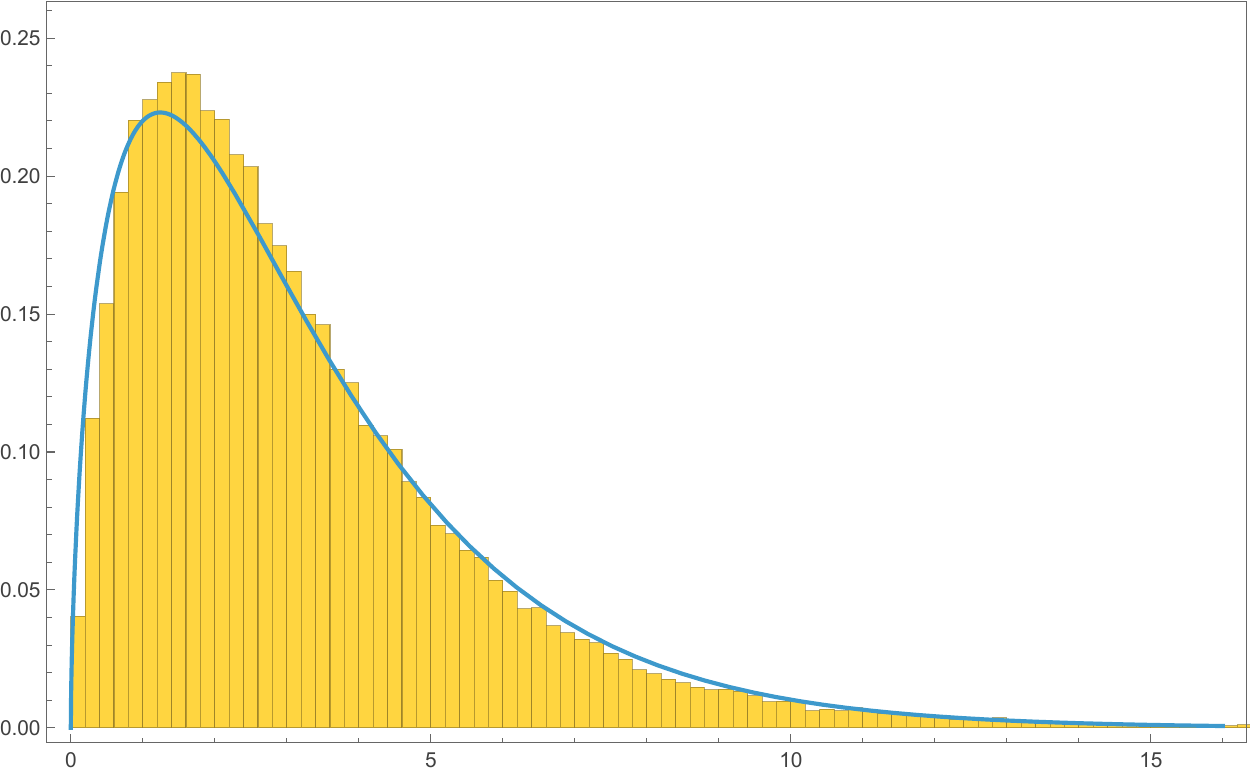}
  \caption{Histogram of $b_{1,2}/c_{50, 2}$ based on $10^5$ replications
at $n=50$, with the density of $\chi^2_{\nu_{50,2}}$ superimposed.}
  \label{fig:momentmatched}
\end{minipage}
\hfill
\begin{minipage}[t]{0.48\textwidth}
  \centering
  \includegraphics[width=\linewidth]{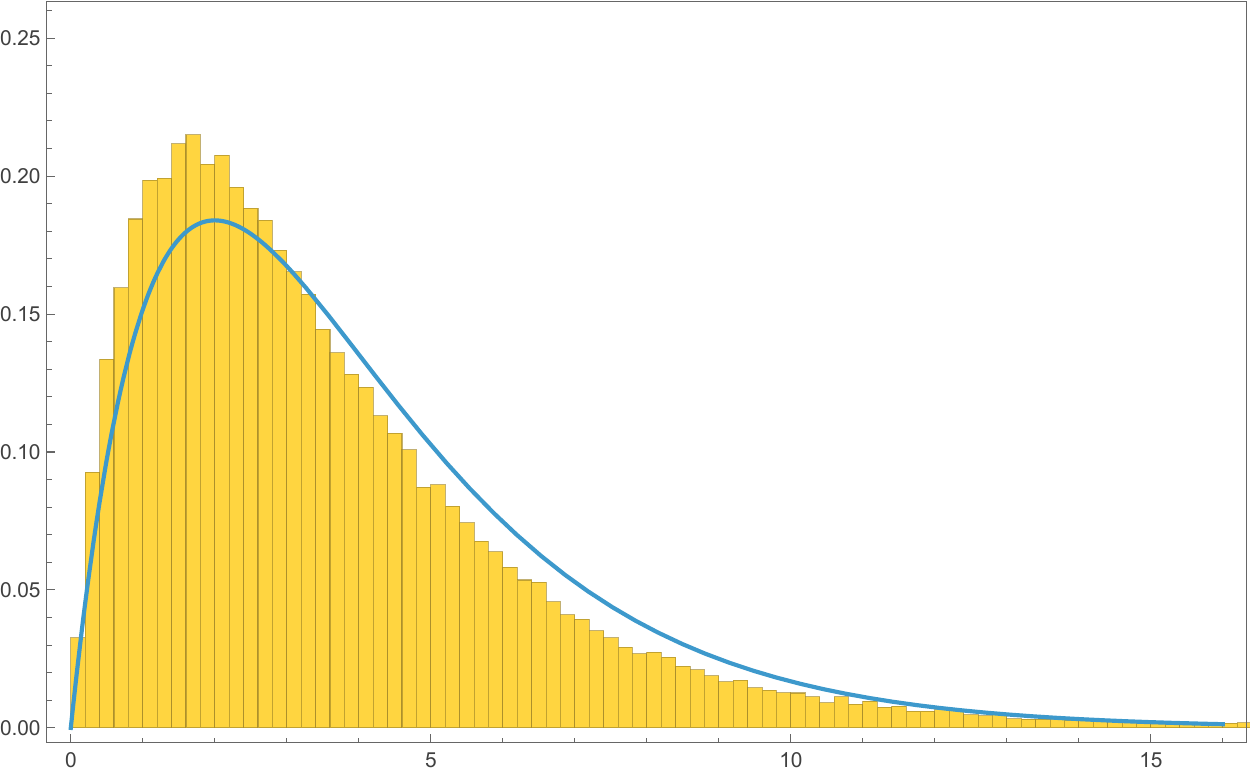}
  \caption{Histogram of $50b_{1,2}/6$ based on $10^5$ replications at $n=50$, with
the density of $\chi^2_4$ superimposed.}
  \label{fig:classical}
\end{minipage}
\end{figure}

For $p=2$, Figures~\ref{fig:momentmatched} and \ref{fig:classical} compare the
corresponding chi-square densities with
simulated histograms based on $10^5$
replications at $n=50$.
Figure~\ref{fig:momentmatched} illustrates the moment-matched
approximation \eqref{eq:mm_chisq}, whereas
Figure~\ref{fig:classical} illustrates the classical approximation
\eqref{eq:mardia_classical}.

Table~\ref{tab:rejection} reports the empirical upper-tail
rejection probabilities for the three
approximations, based on $10^6$ Monte Carlo
replications.
At the nominal level of $0.05$, the moment-matched
approximation gives rejection probabilities
closer to the nominal level than both
alternatives for all sample sizes considered.
At the nominal level of $0.01$, it improves
on the mean-corrected approximation for all
sample sizes considered, but neither the
moment-matched nor the classical approximation
is uniformly superior to the other.

\begin{table}[tbh]
\centering
\caption{Empirical null rejection probabilities for the bivariate
skewness test (based on $10^6$ replications) using approximations 
\eqref{eq:mm_chisq}, \eqref{eq:mardia_classical}, and \eqref{eq:mardia_corrected}.}
\label{tab:rejection}
\begin{tabular}{rcccccc}
\toprule
$p=2$
& \multicolumn{3}{c}{Nominal level $0.05$}
& \multicolumn{3}{c}{Nominal level $0.01$}\\
\cmidrule(lr){2-4}\cmidrule(lr){5-7}
$n$
& \eqref{eq:mm_chisq}
& \eqref{eq:mardia_classical}
& \eqref{eq:mardia_corrected}
& \eqref{eq:mm_chisq}
& \eqref{eq:mardia_classical}
& \eqref{eq:mardia_corrected}
\\
\midrule
20  & 0.0470 & 0.0262 & 0.0604
    & 0.0112 & 0.0061 & 0.0177\\
30  & 0.0467 & 0.0353 & 0.0602
    & 0.0114 & 0.0090 & 0.0179\\
50  & 0.0474 & 0.0424 & 0.0585
    & 0.0114 & 0.0113 & 0.0171\\
100 & 0.0481 & 0.0474 & 0.0562
    & 0.0110 & 0.0117 & 0.0148\\
\bottomrule
\end{tabular}
\end{table}

\begin{figure}[btp]
\centering
\begin{minipage}[t]{0.48\textwidth}
  \centering
  \includegraphics[width=\linewidth]{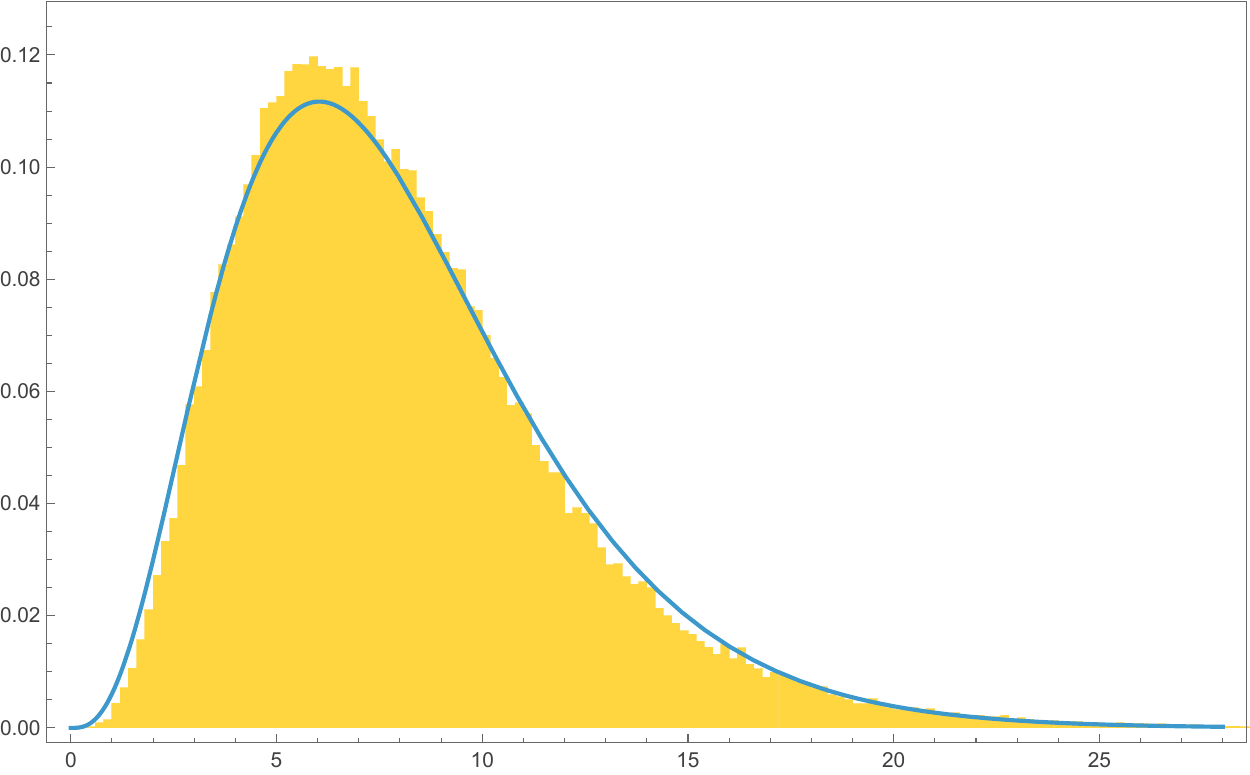}
  \caption{Histogram of $b_{1,3}/c_{50, 3}$ based on $10^5$ replications
at $n=50$, with the density of $\chi^2_{\nu_{50,3}}$ superimposed.}
  \label{fig:momentmatched_p3}
\end{minipage}
\hfill
\begin{minipage}[t]{0.48\textwidth}
  \centering
  \includegraphics[width=\linewidth]{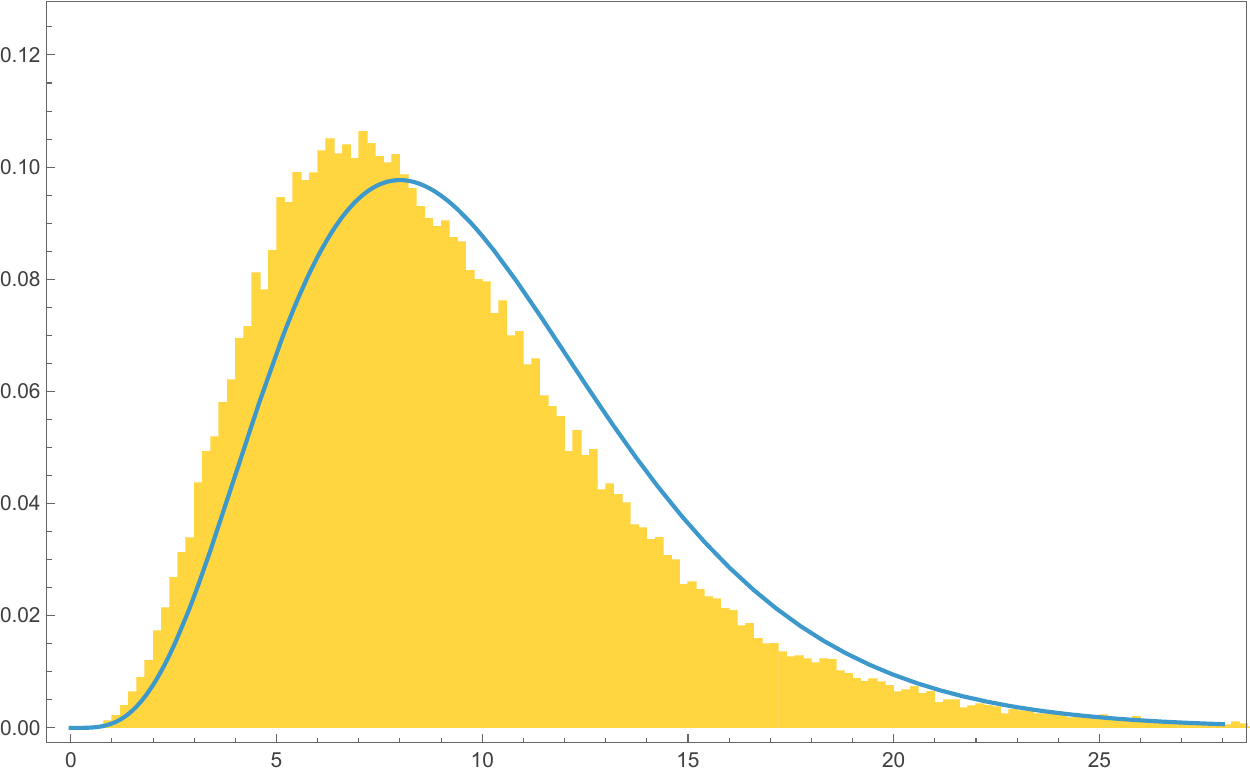}
  \caption{Histogram of $50b_{1,3}/6$ based on $10^5$ replications at $n=50$, with
the density of $\chi^2_{10}$ superimposed.}
  \label{fig:classical_p3}
\end{minipage}
\end{figure}

For \(p=3\), Figures~\ref{fig:momentmatched_p3}
 and \ref{fig:classical_p3} provide the corresponding comparisons at \(n=50\), using \(10^5\) Monte Carlo replications. 
Figure~\ref{fig:momentmatched_p3} shows the histogram of \(b_{1,3}/c_{50,3}\) together with 
the density of \(\chi^2_{\nu_{50,3}}\), whereas Figure~\ref{fig:classical_p3} shows the histogram of \(50b_{1,3}/6\) 
together with the density of \(\chi^2_{10}\). 
The moment-matched approximation appears to follow the simulated distribution more closely in these graphical comparisons.

Table~\ref{tab:trivariate_rejection} reports the empirical upper-tail rejection probabilities for \(p=3\), 
based on \(10^6\) Monte Carlo replications for each sample size. 
At the nominal level of $0.05$, the moment-matched approximation provides the closest rejection 
probability to the nominal level among the three methods for every sample size considered. 
The classical approximation is conservative for small \(n\), whereas the mean-corrected approximation tends to be liberal. 
At the nominal level of $0.01$, the moment-matched approximation improves on the mean-corrected approximation for all four sample sizes, 
but is not uniformly superior to the classical approximation. 
This indicates that matching two moments alone does not guarantee accurate upper-tail probabilities.

\begin{table}[tbh]
\centering
\caption{
Empirical null rejection probabilities for the trivariate
skewness test (based on $10^6$ replications) using approximations 
\eqref{eq:mm_chisq}, \eqref{eq:mardia_classical}, and \eqref{eq:mardia_corrected}.}
\label{tab:trivariate_rejection}
\begin{tabular}{rcccccc}
\hline
$p=3$ & \multicolumn{3}{c}{Nominal level $0.05$} & \multicolumn{3}{c}{Nominal level $0.01$} \\
\cline{2-4}\cline{5-7}
$n$
& \eqref{eq:mm_chisq}
& \eqref{eq:mardia_classical}
& \eqref{eq:mardia_corrected}
& \eqref{eq:mm_chisq}
& \eqref{eq:mardia_classical}
& \eqref{eq:mardia_corrected}
\\
\hline
20 & 0.0502 & 0.0213 & 0.0627 & 0.0125 & 0.0048 & 0.0186 \\
30 & 0.0491 & 0.0324 & 0.0644 & 0.0124 & 0.0085 & 0.0198 \\
50 & 0.0490 & 0.0418 & 0.0634 & 0.0122 & 0.0114 & 0.0191 \\
100 & 0.0497 & 0.0483 & 0.0602 & 0.0115 & 0.0122 & 0.0162 \\
\hline
\end{tabular}
\end{table}

\FloatBarrier
\section{Concluding remarks}\label{sec:conclusion}

Under multivariate normality, polar decomposition shows that Mardia's sample skewness and kurtosis are each independent of the 
sample covariance matrix. 
The resulting Geary-type expectation identities hold in arbitrary dimension 
and yield exact expressions for the variance of sample skewness and the third central moment of sample kurtosis for \(p=2,3\).

The exact mean and variance of sample skewness also yield a moment-matched chi-square approximation to its null distribution. 
For both dimensions and all sample sizes considered, the empirical rejection probabilities are closer to the nominal 
$5\%$ level than those obtained using Mardia's classical and mean-corrected approximations. 
At the $1\%$ level, the approximation improves on the mean-corrected method 
but is not uniformly more accurate than the classical method.

The results for \(p=2,3\), together with the univariate (\(p=1\)) moment recurrence formulas of \cite{nakagawa2016recurrence}, 
may guide future work on conjecturing and proving general formulas for \(\operatorname{Var}[b_{1,p}]\) and \(\mu_3[b_{2,p}]\).

\section*{Disclosure}
During the preparation of this work, the authors used ChatGPT (GPT-5.6 Sol, OpenAI) to improve computer programs 
for deriving moments of Mardia's sample skewness and kurtosis and to explore an alternative proof of Theorems 1 and 2. 
In particular, ChatGPT suggested the use of Basu's theorem as an alternative approach 
to establishing the independence results. 
The authors independently verified the mathematical validity of this argument. 
After using this tool/service, the authors reviewed and edited the content 
as needed and take full responsibility for the content of the publication.

\section*{Acknowledgments}
This work was partially supported by JSPS KAKENHI Grant No. 26K14751.


\begin{thebibliography}{99}

\bibitem[Basu(1955)]{Basu1955}
Basu, D., 1955.
On statistics independent of a complete sufficient statistic.
Sankhy\=a 15, 377--380.

\bibitem[Chikuse(2003)]{Chikuse2003}
Chikuse, Y., 2003.
Statistics on Special Manifolds. Springer, New York.
\url{https://doi.org/10.1007/978-0-387-21540-2}

\bibitem[Geary(1933)]{Geary1933}
Geary, R.C., 1933.
A general expression for the moments of certain symmetrical
functions of normal samples.
Biometrika 25, 184--186.
\url{https://doi.org/10.2307/2332273}

\bibitem[Mardia(1970)]{Mardia1970}
Mardia, K.V., 1970.
Measures of multivariate skewness and kurtosis with applications.
Biometrika 57, 519--530.
\url{https://doi.org/10.1093/biomet/57.3.519}

\bibitem[Mardia(1974)]{Mardia1974}
Mardia, K.V., 1974.
Applications of some measures of multivariate skewness and kurtosis
in testing normality and robustness studies.
Sankhy\=a: The Indian Journal of Statistics, Series B 36, 115--128.

\bibitem[Mardia and Foster(1983)]{MardiaFoster1983}
Mardia, K.V., Foster, K., 1983.
Omnibus tests of multinormality based on skewness and kurtosis.
Communications in Statistics---Theory and Methods 12(2), 207--221.
\url{https://doi.org/10.1080/03610928308828452}

\bibitem[Mardia and Kanazawa(1983)]{MardiaKanazawa1983}
Mardia, K.V., Kanazawa, M., 1983.
The null distribution of multivariate kurtosis.
Communications in Statistics---Simulation and Computation
12(5), 569--576.
\url{https://doi.org/10.1080/03610918308812343}

\bibitem[Muirhead(1982)]{Muirhead1982}
Muirhead, R.J., 1982.
Aspects of Multivariate Statistical Theory. Wiley, New York.
\url{https://doi.org/10.1002/9780470316559}

\bibitem[Nakagawa et al.(1998)]{NakagawaNikiHashiguchi1998}
Nakagawa, S., Niki, N., Hashiguchi, H., 1998.
Computer algebra application to the distribution of
sample correlation coefficient.
Mathematics and Computers in Simulation 45(1--2), 23--32.
\url{https://doi.org/10.1016/S0378-4754(97)00083-9}

\bibitem[Nakagawa et al.(2016)]{nakagawa2016recurrence}
Nakagawa, S., Hashiguchi, H., Niki, N., 2016.
Recurrence formula on the joint distribution of sample skewness
and kurtosis under normality.
Communications in Statistics---Theory and Methods 45(4), 1148--1155.
\url{https://doi.org/10.1080/03610926.2013.859704}

\end{thebibliography}
\end{document}